\documentclass{amsart}
\usepackage{amsmath,amsthm,amssymb,amsfonts}

\makeatletter

\def\proof{\@ifstar{P\,r\,o\,o\,f}{P\,r\,o\,o\,f.\ }}
\renewcommand\th@remark {%
  \thm@headfont{\bfseries}%
  \normalfont % body font
  \thm@preskip\topsep \divide\thm@preskip\tw@
  \thm@postskip\thm@preskip
}

\renewenvironment{equation}{\refstepcounter{equation}$$}{\eqno{(\thesection.\theequation)}$$}
\makeatother

\newcounter{Th}[section]  \newcounter{Lm}[section]
\newcounter{Remark}[section] \newcounter{Example}[section] \newcounter{Def}[section]
\newcounter{lcounter}[section]
\newcounter{Corol}[section]

\newenvironment{Th}[1][\relax]
    {\medspace\refstepcounter{Example}T\,h\,e\,o\,r\,e\,m \arabic{section}.\theExample.\ \it}
    {\rm\medspace}

\newenvironment{Th.}[1][\relax]
    {\medspace\refstepcounter{Example}T\,h\,e\,o\,r\,e\,m \arabic{section}.\theExample.\ \it}
    {\rm\medspace}

\newenvironment{Pr.}[1][\relax]
    {\medspace\refstepcounter{Example}P\,r\,o\,p\,o\,s\,i\,t\,i\,o\,n \arabic{section}.\theExample.\ \it}
    {\rm\medspace}

\newenvironment{Lm}[1][\relax]
    {\medspace\refstepcounter{Example}L\,e\,m\,m\,a \arabic{section}.\theExample.\ \it}
    {\rm\medspace}

\newenvironment{Corol}[1][\relax]
    {\medspace\refstepcounter{Example} C\,o\,r\,o\,l\,l\,a\,r\,y \arabic{section}.\theExample.\ \it}
    {\rm\medspace}

\newenvironment{Remark}[1][\relax]
    {\medspace\refstepcounter{Example}R\,e\,m\,a\,r\,k \arabic{section}.\theExample.\rm\ }
    {\medspace}
\newenvironment{Example}[1][\relax]
    {\medspace\refstepcounter{Example}E\,x\,a\,m\,p\,l\,e \arabic{section}.\theExample.\rm\ }
    {\medspace}

\newenvironment{Def}[1][\relax]
    {\medspace\refstepcounter{Example}D\,e\,f\,i\,n\,i\,t\,i\,o\,n \arabic{section}.\theExample.\rm\ }
    {\medspace}

\newenvironment{Def.}[1][\relax]
    {\medspace\refstepcounter{Example}D\,e\,f\,i\,n\,i\,t\,i\,o\,n \arabic{section}.\theExample.\rm\ }
    {\medspace}

\def\au#1{\emph{#1}}% Написание авторов в библиографии
\def\tit#1{{#1}}% Написание заголовков в библиографии
\def\R{{\mathbb R}}  % множество вещественных чисел
\def\cl{\mathop{\rm cl\,}}

\def\VolE#1{{\bf #1}}

\def\co{\mbox{\rm co}\,}
\def\l{\langle}
\def\r{\rangle}
\def\ll{\lambda}
\def\d {\partial\,}
\def\ep{\varepsilon} % эпсилон
\def\Int {\mbox{\rm int\,}} % внутренность

\def\L {{\mathcal L}}  % каллиграфическое L
\def\diam {\mbox{\rm diam}\,}
\def\diff {\,\frac{\,*\,}{}\,}

\begin{document}

\title[On the splitting problem for selections]{On the splitting problem for selections}

\author[M. V. Balashov and D. Repov\v{s}]{Maxim V. Balashov and Du\v{s}an Repov\v{s}}

%\thanks{Supported by }

\address{Department of Higher Mathematics, Moscow Institute of Physics and Technology, Institutski str. 9,
Dolgoprudny, Moscow region, Russia 141700}
\email{balashov@mail.mipt.ru}
\address{Faculty of Mathematics and
Physics, and Faculty of Education,
University of Ljubljana, P.O.B. 2964, Ljubljana, Slovenia 1001}
\email{dusan.repovs@guest.arnes.si}



\begin{abstract}
We investigate  when does the Repov\v{s}-Semenov Splitting 
problem for selections
have an
affirmative solution for continuous
set-valued mappings in finite-dimensional Banach spaces. We prove
that this happens
when images of set-valued mappings or even
their graphs  are P-sets (in the sense of Balashov) 
or strictly convex sets. We
also consider
an
example which shows that there is no affirmative
solution of this problem even in the simplest case 
in $\R^{3}$. We also obtain affirmative
solution of the Approximate splitting problem for Lipschitz
continuous selections in the
Hilbert space.
\end{abstract}

\date{\today}
\maketitle

\section{Introduction}

The splitting problem for selections was recently 
stated in
\cite{Repovs+Semenov}. Let $F_{i}:X\to 2^{Y_{i}}$, $i=1,2$, be
any (lower semi)
continuous mappings with closed convex images and
let $L:Y_{1}\oplus Y_{2}\to Y$ be any linear surjection. The
splitting problem is the problem of representation of any
continuous selection $f(x)\in L(F_{1}(x),F_{2}(x))$ in the form
$f(x)= L(f_{1}(x),f_{2}(x))$, where $f_{i}(x)\in F_{i}(x)$ are
some continuous selections, $i=1,2$.

This problem is related to some classical problems of set-valued
analysis. First, it is a special case of the
selection problem which
is sufficiently common for various applications
\cite{RepovsSemenovBook}. In particular it is quite close to the
celebrated
problem of parametrization of set-valued mappings
\cite{Aubin+Frankowska,Ivanov+Balashov,Ornelas}.

Second, every affirmative solution of this problem is in fact, an
answer to the following
question: When does the operation of
intersection of 
two (continuous) set-valued mappings yield
a
continuous (or lower semicontinuous) set-valued mapping with
respect to the Hausdorff metric? This question is also quite
common for certain
branches of set-valued and nonsmooth analysis.

It is well known that the
intersection of two continuous
set-valued mappings is not necessarily continuous
\cite{Aubin+Frankowska}. We shall first
consider the
extreme example which demonstrates the last assertion.

Consider Question 4.6 from \cite{Repovs+Semenov}: Do there exist for
every closed convex sets $A$, $B$ and $C=A+B$, continuous
functions $a:C\to A$ and $b:C\to B$ with the property that
$a(c)+b(c)=c$, for all $c\in C$?  In a space of
dimension $\ge 3$ the answer is negative.

\begin{Example}\label{Ex1}
Consider the following sets in the $3$-dimension Euclidean space $\R^{3}$ 
(where $\co(X)$ denotes the convex hull of $X$):
$$
D_{0}=\{ (\cos t,\sin t, 0)\ |\ t\in [0,\pi]\},\quad
A_{0}=\co(D_{0}\cup\{ (1,0,1)\}\cup\{ (-1,0,1)\}),
$$
$$
D_{1}=\{ (\cos t,\sin t, 1)\ |\ t\in [-\pi,0]\},\quad
A_{1}=\co(D_{1}\cup\{ (1,0,0)\}\cup\{ (-1,0,0)\}),
$$
and $A=A_{0}\cup A_{1}$. It's easy to see that $A$ is a convex
compact set. We also define the set $B=\co((0,0,0),(0,0,1))$ and
the set $C=A+B$.

Let $\Gamma =\{ (\cos t,\sin t, 1-\frac{2}{\pi}t)\ |\ t\in
[-\frac{\pi}{2},\frac{\pi}{2}]\}$,
$\Gamma \subset \d C$. Let
$$
\Gamma_{1} =\left\{ \left(\cos t,\sin t, 1-\frac{2}{\pi}t\right)\
\left| \ t\in (0,\pi/2]\right.\right\},
$$
$$
\Gamma_{2} =\left\{ \left(\cos t,\sin t, 1-\frac{2}{\pi}t\right)\
\left| \ t\in [-\pi/2,0)\right.\right\}
$$
and $c_{0}=(1,0,1)$.

Suppose that $c\in \Gamma_{1}$. In this case there exists
only one
pair of points $a(c)\in A$ and $b(c)\in B$ with the property
$a(c)+b(c)=c$. Indeed, if $t_{c}\in (0,\frac{\pi}{2}]$ satisfies
$c=(\cos t_{c},\sin t_{c}, 1-\frac{2}{\pi}t_{c})$ then
$a(c)=(\cos t_{c}, \sin t_{c},0)$, $b(c)=(0,0,
1-\frac{2}{\pi}t_{c})$. The point $a(c)$ is unique because it is
an exposed point of the set $A$ for vector $p_{c}=(\cos t_{c},\sin
t_{c}, 0)$. Clearly, the point $b(c)$ is also unique. So we have
\begin{equation}
 \lim\limits_{c\to c_{0},\ c\in
\Gamma_{1}}a(c)=(1,0,0).\label{1c}
\end{equation}

In the case when $c\in \Gamma_{2}$, similar considerations show
that there
exists only one pair of points $a(c)=(\cos t_{c}, \sin
t_{c},1)\in A$ and $b(c)=(0,0,-\frac{2}{\pi}t_{c})\in B$ with
$a(c)+b(c)=c$ and such that

\begin{equation}
\lim\limits_{c\to c_{0},\ c\in \Gamma_{2}}a(c)=(1,0,1).\label{2c}
\end{equation}

Formulae
(1.\ref{1c}) and (1.\ref{2c}) 
show that $a(c)$ is not
continuous at the point $c=c_{0}$.

\end{Example}

Simultaneously, we want to emphasize that the set-valued mappings
\begin{equation}\label{bad-1}
C\ni c\to (c-B)\cap A
\end{equation}
and
\begin{equation}\label{bad-2}
C\ni c\to (A,B)\cap L^{-1}(c)
\end{equation}
 do not allow any
 continuous (on $c\in C$) selection. Here $L^{-1}(c)=\{ (y_{1},y_{2})\in \R^{3}\times \R^{3}\ |\
 y_{1}+y_{2}=c\}$. 
 
 Indeed, otherwise in the case (1.\ref{bad-1}) we could
 choose this selection as  $a(c)\in
 (c-B)\cap A\subset A$ and set
 $b(c)=c-a(c)\in B$. In the case (1.\ref{bad-2}) we could choose $(a(c),b(c))\in
(A,B)\cap L^{-1}(c)$. In both cases we would have $a(c)+b(c)=c$,
for all
$c\in C$. This would contradict
the fact that function $a(c)$, as it
follows by Example 1.1, is not continuous at the point
$c=c_{0}=(1,0,1)$.

We shall
further obtain an
affirmative
solution of the splitting
problem for selections for some special cases in
finite-dimensional Banach spaces. It suffices to solve this
problem in Euclidean space $\R^{n}$ (with the inner product
$\l\cdot,\cdot\r$) because all norms in any finite-dimensional
Banach space are equivalent.

Our general idea is to  prove continuity of the intersection
$(F_{1}(x),F_{2}(x))\cap L^{-1}(f(x))$. This is the key idea. When
this is done we can choose some continuous selection (e.g. the
Steiner
point) of the map $x\to (F_{1}(x),F_{2}(x))\cap L^{-1}(f(x))$ and
solve the problem.

Unfortunately, this map is not continuous even in the simplest
cases (as we can see from Example 1.1). So we need to invoke some
special geometrical properties of maps $F_{i}$ or surjection $L$.

The main results of Section 2 are Theorem~2.\ref{affsection} with Corollary~2.\ref{sect}, Theorem~2.\ref{Ex3} and Theorem~2.\ref{srtictlyConv}. The key geometric objects
used in Section 2 are P-sets (see Definition~2.\ref{p-set}) 
and strictly convex sets.

In Section 3 we shall consider situation when in
infinite-dimensional Hilbert spaces $X$, $Y_{1}$, $Y_{2}$ there exist
for
every $\ep>0$ and any Lipschitz continuous selection $f(x)\in
L(F_{1}(x),F_{2}(x))$,  Lipschitz continuous selections
$f_{i}(x)\in F_{i}(x)+B^{Y_{i}}_{\ep}(0)$, $i=1,2$, with the
property
$f(x)=L(f_{1}(x),f_{2}(x))$, for all $x$. Here,
$B_{\ep}^{Y_{i}}(0)=\{ y\in Y_{i}\ |\ \| y\|\le \ep\}$. The main
results of Section 3 are Theorem~3.\ref{unicont} and
Theorem~3.\ref{ApproxSolution}.

We need to give some definitions for further explanation.
We shall
say that the
subspace $\L\subset \L_{1}\oplus\L_{2}$ is {\it
not parallel} to the subspaces $\L_{1}$ and $\L_{2}$ if for any
pair of distinct
points $w_{1},w_{2}\in \L$,
the
projections of $w_{1}$ and
$w_{2}$ onto 
$\L_{1}$ (resp. $\L_{2}$) parallel to
$\L_{2}$ (resp. $\L_{1}$)
yield different points. 

Let $h$ be the {\it Hausdorff distance}. For any bounded subsets
$A,B$ of a Banach space $X$ we have
$$
h(A,B)=\inf\{ h>0\ |\ A\subset B+B^{X}_{h}(0),\quad B\subset
A+B^{X}_{h}(0)\}.
$$

For any subsets $A,B$ of a linear space $X$ the operation
$$
A\diff B=\{x\in X\ |\ x+B\subset A\} = \bigcap\limits_{b\in
B}(A-b)
$$
is called the {\it geometric difference} (or the 
{\it Minkowski-Pontryagin difference}) of sets $A$ and $B$. A direct consequence of
the
definition of geometric difference is that $(A\diff B)+B\subset
A$.

The {\it Chebyshev center} $c(A)$ of a
convex closed bounded subset $A$
of a
Hilbert space $X$ is the point
$$
c(A) = \arg\inf\limits_{x\in X}\left(\sup\limits_{a\in A}\|
x-a\|\right).
$$
Chebyshev center always exists and it is unique.

Let $X,Y$ be any
Banach spaces. We say that a continuous linear
surjection $L:X\to Y$ has the
{\it Lip-property} if the set-valued
mapping $L^{-1}(y)=\{ x\in X\ |\ Lx=y\}$ has a
Lipschitz (at $y$)
selection $l(y)\in L^{-1}(y)$.

For example, if $Y_{1}=Y_{2}=Y$ and $L:Y_{1}\oplus Y_{2}\to Y$,
$L(y_{1},y_{2})=y_{1}+y_{2}$, then $l(y)=\left(\frac12 y,\frac12
y\right)$.

\bigskip

\section{P-sets and the splitting problem for selections}

We shall
obtain
an affirmative
solution in certain
cases when the
images of the
set-valued mapping are P-sets (\cite{Bal}). 
Let $q\in\R^{n}$ be an
arbitrary unit vector and $L(q)=\{ z\in\R^{n}\ |\ \l q,z\r=0\}$,
$l(q)=\{ \ll q\ |\ \ll\in\R\}$. The space $\R^{n}$ is the
orthogonal sum of sets $L(q)$ and $l(q)$: $\R^{n}=L(q)\oplus
l(q)$. Any point $z\in\R^{n}$ can be expressed in the form
$z=x+\mu q$, where $\mu\in \R$, $x\in L(q)$, or $z=(x,\mu)$. Let
$P_{q}:\R^{n}\to L(q)$ be the operator of orthogonal projection: for
any $z\in\R^{n}$, $P_{q}z=x$, where $z=(x,\mu)$.

Let $A\subset\R^{n}$ be any convex compact set. Let's define the
function $f_{A,q}:P_{q}A\to \R$ by
\begin{equation}\label{function}
f_{A,q}(x)=\min\{ \mu\ |\ (x,\mu )\in A\},\qquad\mbox{for all}\ 
x\in P_{q}A.
\end{equation}

\begin{Def}
  \label{p-set} (\cite{Bal})
A convex compact subset $A\subset\R^{n}$ is called a
 \it $P$-set\rm, if for any unit vector  $q$, the function
  $f_{A,q}$ (2.\ref{function}) is continuous on the set $P_{q}A$.
\end{Def}

\begin{Pr.} 
\label{p-set-properties} 
(\cite{Bal,Polovinkin&Balashov})
\label{ThisWas}
Any convex compact subset of $\R^{2}$ is a $P$-set. In $\R^{n}$
any convex polyhedron is a $P$-set, any strictly convex compact subset
is a $P$-set, and any finite Minkowski sum of $P$-sets is a $P$-set.
If $L:\R^{n}\to\R^{m}$ is a linear operator and $A\subset \R^{n}$
is a $P$-set then $LA\subset \R^{m}$ is also a $P$-set. Moreover,
the map $L:A\to LA$ is  open  in induced topologies.
\end{Pr.}

  We emphasize
that a $P$-set is not necessarily a polyhedron or strictly convex.

\begin{Example}\label{Ex2}
The cylinder $\{ (x_{1},x_{2},x_{3})\ |\ x_{1}^{2}+x_{2}^{2}\le
1,\ 0\le x_{3}\le 1 \}$ is a P-set as the Minkowski sum of two P-sets
$$
\{ (x_{1},x_{2},0)\ |\ x_{1}^{2}+x_{2}^{2}\le 1\}+\{ (0,0,x_{3})\
|\ 0\le x_{3}\le 1 \}.
$$

On the other hand, the subset $A_{0}\subset\R^{3}$,
$$
A_{0}=\co\left( \{ (x_{1}-1)^{2}+x_{2}^{2}=1,\ x_{3}=1\}\cup \{
(0,0,0)\} \right),
$$
is not a $P$-set.

Indeed, for  $q=(0,0,1)$ it is easy to see that $f_{A_{0},q}$ is
not upper semicontinuous at the point $(0,0)\in P_{q}A_{0}$. Note
that the sum $A_{0}+A_{1}$ (where $A_{1}$ is an arbitrary convex compact set)
is not a $P$-set \cite{Bal}.
\end{Example}

It was proved
in
\cite{Bal}
that if the subset
$A\subset\R^{n}$ is convex and compact then the function $f_{A,q}$
(2.\ref{function}) is lower semicontinuous on $P_{q}A$. This is
quite obvious. Therefore
the question about continuity of the function
$f_{A,q}$ is the question about its
upper semicontinuity.

\def\dom{\mbox{\rm dom}\,}

The {\it domain} of the set-valued mapping $F:\R^{n}\to 2^{\R^{m}}$ is the
set $\dom F=\{ x\in \R^{n}\ |\ F(x)\ne\emptyset\}$.

\begin{Th}
\label{affsection} Let $A\subset \R^{n}$ be any convex compact
subset and $\L\subset \R^{n}$
any
subspace. Suppose that one of the following properties is satisfied:\\
 1) $\dim \L=n-1$, or \\
 2) $A$ is a $P$-set.\\
 Then the set-valued mapping $F(z)=(z+\L)\cap A$ is continuous in the Hausdorff metric.% on $\dom F$.
\end{Th}

Remark: It is easy to see that $\dom F = A+\L$.

\proof 1) Part (1) is a well known fact, but we prove it for
completeness.
Let $A\subset \R^{n}$ be an arbitrary convex compact subset and
$\dim \L = n-1$. Consider $F(z)=(z+\L)\cap A$, for all $z\in \dom
F$. Using the Closed graph theorem for set-valued mappings
(\cite[Theorem~8.3.1 ]{Aubin}) we conclude
that map $F$ is
upper
semicontinuous at any point $z_{0}\in \dom F$. 

This means that for
any sequence  $\{ z_{k}\}\subset \dom F$, $z_{k}\to z_{0}$, and
any $\ep>0$,
there exists a number $k_{\ep}$ such that for any
$k>k_{\ep}$ the following holds:
\begin{equation}\label{upper}
 F(z_{k})\subset F(z_{0})+B_{\ep}(0).
\end{equation}
Suppose that lower semicontinuity fails at some point $z_{0}\in
\dom F$. Then there exist
a number $\ep_{0}>0$ and a sequence
$\{ z_{k}\}\subset \dom F$
such
that $\lim z_{k}=z_{0}$ and
\begin{equation}\label{lower}
 F(z_{0})\not\subset F(z_{k})+B_{\ep_{0}}(0).
\end{equation}
This implies that $z_{k}\notin z_{0}+\L$, for all $k$. From the
condition (2.\ref{lower}) we conclude
that for any $k$, there exists a
point $w_{k}\in F(z_{0})$ for which $w_{k}\notin
F(z_{k})+B_{\ep_{0}}(0)$. 

We may
assume that $w_{k}\to w_{0}\in
F(z_{0})$,
due to the compactness of the set $F(z_{0})$ and since
\begin{equation}\label{lower1}
 w_{0}\notin F(z_{k})+B_{\ep_{0}/2}(0).
\end{equation}
Without loss of generality we may
assume that $z_{1}\in A$.
Otherwise, we can choose an arbitrary point $\tilde z_{1}\in
F(z_{1})$ instead of $z_{1}$. 
We can
also suppose that $\|
z_{1}-w_{0}\|>\frac{\ep_{0}}{2}$, otherwise we could reduce
$\ep_{0}$. 

Let $\varphi$ be the angle between the segment
$[w_{0},z_{1}]$ and the hyperplane $\L$, $\varphi \in
(0,\frac{\pi}{2}]$ (note that the segment $[w_{0},z_{1}]$ is not
parallel to $\L$). Without loss of generality we may assume
that
the halfspace with the bound $z_{0}+\L$, which contains $z_{1}$,
contains the entire sequence $\{ z_{k}\}$.

Let us choose a number $k$ for which the distance from the point
$z_{k}$ to the hyperplane $z_{0}+\L$ is sufficiently small:
$$
\varrho (z_{k},z_{0}+\L)<\frac{\ep_{0}}{2}\sin\varphi.
$$
Define the point $w=w_{0}+\frac{z_{1}-w_{0}}{\|
z_{1}-w_{0}\|}\frac{\ep_{0}}{2}$. We have
$$
\varrho (w,z_{0}+\L)=\frac{\ep_{0}}{2}\sin\varphi > \varrho
(z_{k},z_{0}+\L).
$$
We can now conclude that the points $w_{0}$ and $w$ lie on the
opposite sides of the hyperplane $z_{k}+\L$. This follows from
the last estimate, the
inclusion $w_{0}\in z_{0}+\L$ and the fact
that the
points $z_{k}$, $z_{1}$ (and consequently $w$) lie on the
same side of the hyperplane $z_{0}+\L$. 
This means that the following holds:
\begin{equation}\label{Gr}
 [w_{0},w]\cap (z_{k}+\L)\ne\emptyset.
\end{equation}
From the inclusions $w_{0}\in F(z_{0})$, $z_{1}\in A$ we obtain
that $[w_{0},w]\subset [w_{0},z_{1}]\subset A$. We can conclude
from this inclusion and the equality $\|
w_{0}-w\|=\frac{\ep_{0}}{2}$ that $[w_{0},w]\subset A\cap
B_{\ep_{0}/2}(w_{0})$. According to (2.\ref{Gr}), we have
$$
 A\cap B_{\ep_{0}/2}(w_{0})\cap (z_{k}+\L)\ne\emptyset,
$$
i.e. $w_{0}\in \left( A\cap (z_{k}+\L)\right)+B_{\ep_{0}/2}(0)$.
This contradicts the existence of the inclusion (2.\ref{lower1}).

The upper and lower semicontinuity imply continuity in the
Hausdorff metric.

(2) Suppose now that $\dim \L=m$, $1\le m\le n$, and that
$A\subset
\R^{n}$ is a $P$-set. Upper semicontinuity can be proved in the
same way as in (1) above.

Assuming to the contrary, as in (1), we conclude that there exist
$\ep_{0}>0$, $w_{0}\in F(z_{0})$ and a sequence $\{ z_{k}\}\subset
\dom F$, $\lim z_{k}=z_{0}$ such that
$$
w_{0}\notin F(z_{k})+B_{\ep_{0}/2}(0).
$$
The map $F$ is upper semicontinuous and not continuous. Hence
$$
F_{0}=\lim\sup\limits_{k\to\infty}F(z_{k})=\bigcap\limits_{\ep>0}\bigcap\limits_{\delta>0}
\bigcup\limits_{k>\delta}\left( F(z_{k})+B_{\ep}(0)\right)\subset
F(z_{0}).
$$
Let us fix any point $w\in F_{0}$. Obviously, $w\ne w_{0}$.
Let $q=\frac{w-w_{0}}{\| w-w_{0}\|}$. Suppose
that $v_{k_{n}}\in
F(z_{k_{n}})$ is a sequence such that $v_{k_{n}}\to w$. We have
$F(z_{k_{n}})=F(v_{k_{n}})$, $F(z_{0})=F(w_{0})=F(w)$.

Let $w_{0}=(x_{0},\mu_{0})$. Note that $x_{0}=P_{q}w$. Let
$x_{k_{n}}=P_{q}v_{k_{n}}$. For sufficiently large $n$ (when $\|
w-v_{k_{n}} \|<\frac{\ep_{0}}{4}$) we have $(v_{k_{n}}+l(q))\cap
B_{\ep_{0}/4}(w_{0})\ne\emptyset$, $l(q)=\{\ll (w-w_{0})\ |\
\ll\in\R\}\subset \L$
 and consequently, $f_{A,q}(x_{k_{n}})\ge
\mu_{0}+\frac{\ep_{0}}{4}$. Together with $f_{A,q}(x_{0})\le
\mu_{0}$, this contradicts the definition of a $P$-set. \qed

\def\dom{\mbox{\rm dom}\,}
\def\graph{\mbox{\rm graph}\,}
The graph
$\graph F$ of the set-valued mapping $F:\R^{n}\to 2^{\R^{m}}$ is
the set $\graph F = \{ (x,y)\in\R^{n}\times\R^{m}\ |\ x\in\dom F,
\ y\in F(x)\}$.

\begin{Corol} \label{sect}
Let $F:\R^{n}\to 2^{\R^{m}}$ be any
set-valued mapping with a
convex closed graph. \\
(1) If $n=1$ and $\graph F$ is convex and compact then $F$ is
continuous.\\ % on $\dom F$.\\
(2) If $\graph F$ is a $P$-set then $F$ is continuous. % on $\domF$.
\end{Corol}

\proof* of Corollary 2.\ref{sect} follows from
Theorem~2.\ref{affsection} and the
equality
$$
\{ x_{0}\}\times F(x_{0})=\graph F\cap \{
(x,y)\in\R^{n}\times\R^{m}\ |\ x=x_{0}\}.
$$
\qed

Corollary 2.\ref{sect} is false
if $\graph F$ is not a $P$-set.
Let $q$ be the vector $(0,0,1)\in\R^{3}$. Let $\graph F=A_{0}$,
$A_{0}$ be as in
Example 2.\ref{Ex2}, $\dom F=\{
(x_{1}-1)^{2}+x_{2}^{2}\le 1\}$ and $F(x_{1},x_{2})=\{ \lambda \
|\ (x_{1},x_{2},\lambda )\in A_{0}\}$. Then the set-valued
mapping $F$ is not lower semicontinuous at the point $(0,0)$.

The {\it Steiner point}
of  a convex compact subset $A\subset\R^{n}$ is the point
$$
s(A)=\frac{1}{\mu_{n}B_{1}(0)}\int_{\| p\|=1}s(p,A)p\, dp,
$$
$s(p,A)=\sup\limits_{x\in A}\l p,x\r$, where $\mu_{n}$ is the
Lebesgue
measure in $\R^{n}$. For any convex compact subsets $A,B\subset\R^{n}$, we have:
$$
\| s(A)-s(B)\|\le L_{n}h(A,B), \qquad L_{n}=
\frac{2}{\sqrt{\pi}}\frac{\Gamma\left( \frac{n}{2}+1
\right)}{\Gamma\left(\frac{n+1}{2}\right)},
$$
and $s(A)\in A$. The Lipschitz constant $L_{n}$ above is the best
possible \cite{Polovinkin&Balashov,Aubin+Frankowska}.

\begin{Th}\label{Ex3}
Let $A$, $B$ be any closed convex subsets and let $C=A+B$. If $C$ is
a $P$-set then there exist continuous functions $a:C\to A$ and
$b:C\to B$ with the property
that $a(c)+b(c)=c$, for all $c\in C$.
\end{Th}

\proof  Let $L:\R^{n}\times\R^{n}\to \R^{n}$ be a linear
operator, $L(x_{1},x_{2})=x_{1}+x_{2}$. Then $L$ is a surjection. The
set-valued mapping $\R^{n}\ni c\to L^{-1}(c)=\{ (x_{1},x_{2})\ |\
x_{1}+x_{2}=c\}$ is Lipschitz continuous in the Hausdorff metric
\cite[Corollary 3.3.6]{Aubin} and its values are parallel affine
planes of the same dimension. The set $C$ is a $P$-set and
consequently, $A$ and $B$ are $P$-sets, too \cite{Bal}. Using
Corollary 2.\ref{sect}, we conclude
that the set-valued mapping
$C\ni c\to (A,B)\cap L^{-1}(c)$ is continuous. Taking the Steiner
point $s(\cdot)$ of the latter set-valued mapping, we get the following:
$$
(a(c),b(c))=s\left( (A,B)\cap L^{-1}(c)\right).
$$\qed

\begin{Th}\label{Ex4}
 Consider any
  set-valued mappings $F_{1}:\R^{n}\to
2^{\R^{m_{1}}}$ and $F_{2}:\R^{n}\to 2^{\R^{m_{2}}}$. Suppose
that $\graph (F_{1},F_{2})$ is a $P$-set. Suppose that
$L:\R^{m_{1}}\times\R^{m_{2}}\to \R^{k}$ is a linear surjection.
Then for any continuous selection $f(x)\in
L(F_{1}(x),F_{2}(x))$, there exist continuous selections
$f_{i}(x)\in F_{i}(x)$, $i=1,2$, with $f(x)=L(f_{1}(x),f_{2}(x))$.
\end{Th}

\proof We can take
$$
(x,f_{1}(x),f_{2}(x))=s\left(\left( \mbox{graph}\, (F_{1},
F_{2})\right)\cap (x,L^{-1}(f(x)))\right).
$$
The map $x\to (x,L^{-1}(f(x)))$ is continuous, due to
 \cite[Corollary 3.3.6]{Aubin} and the intersection is continuous
by of Corollary 2.\ref{sect}.\qed

\begin{Th}\label{Ex5}
 Suppose that a compact subset $A\subset\R^{n} $ is strictly convex
 and that
$B\subset\R^{n}$ is an arbitrary convex closed subset. In this case
there exist continuous functions $a:C\to A$ and $b:C\to B$ with
the property that
$a(c)+b(c)=c$,
for all $c\in C$.
\end{Th}

\proof We shall
consider $H(c)=(c-A)\cap B$, $c\in C$. Note that
$H(c)\ne \emptyset$, for all $c\in C$. If $\mbox{ri}\,
H(c_{0})\ne\emptyset$, for some $c_{0}\in C$, the condition of
nonempty interior yields the continuity of $H(c)$ at the point
$c=c_{0}$ (\cite{Polovinkin&Balashov,Aubin+Frankowska}). 

Note
that $\mbox{ri}\, H(c_{0})$ is the relative interior of the set
$H(c_{0})$, i.e. the interior of the set $H(c_{0})$ in the
affine hull
of the set $B$.

If $\mbox{ri}\, H(c_{0})=\emptyset$ then $H(c_{0})$ is a single
point, due to the strict convexity of $A$. The intersection of
$c-A$ and $B$ is upper semicontinuous, by the Closed
graph theorem (\cite[Theorem~3.1.8 ]{Aubin}), i.e.
$$
H(c)\subset H(c_{0})+B_{\ep}(0),\qquad\forall c\in
B_{\delta}(c_{0})\cap C.
$$
But this implies
$$
H(c_{0})\in H(c)+B_{\ep}(0),\qquad\forall c\in
B_{\delta}(c_{0})\cap C.
$$
Both formulae yield the continuity (in the Hausdorff metric) at the point
$c=c_{0}$.
Thus $H(c)$ is continuous at any point
$c\in C$ and $b(c)=s(H(c))$,
$a(c)=c-b(c)$.\qed

\begin{Th}\label{Ex6} Let $X$ be any metric space.
Consider any
set-valued
mappings $F_{1}:X\to 2^{\R^{m_{1}}}$ and $F_{2}:X\to
2^{\R^{m_{2}}}$
which are continuous in the Hausdorff metric. 
Suppose that $F_{1}$ has strictly convex compact
images and that
$F_{2}$ has closed convex images. Suppose also
that
$f(x)\in F_{1}(x)+F_{2}(x)$ is a continuous selection. Then there
exist continuous selections $f_{i}(x)\in F_{i}(x)$ with
$f(x)=f_{1}(x)+f_{2}(x)$, for all $x\in X$.
\end{Th}

\proof* is similar to that of Theorem 2.\ref{Ex5}.
The set-valued mapping
$$
H(x)=\left( f(x)-F_{1}(x)\right)\cap F_{2}(x)
$$
is continuous and $f_{2}(x)=s(H(x))$,
$f_{1}(x)=f(x)-f_{2}(x)$.\qed

\begin{Th}\label{srtictlyConv} Let $X$ be any metric space.
Let $F_{i}:X\to2^{\R^{m_{i}}}$, $i=1,2$, be set-valued mappings with strictly convex compact
or single-point images, which is continuous in the
Hausdorff metric.
Let $\L_{i}=\R^{m_{i}}$, $i=1,2$ and let
$L:\L_{1}\oplus\L_{2}\to \R^{k}$ be a linear surjection such that
$\L=\ker L$ is not parallel to $\L_{i}$, $i=1,2$. Then for any
continuous selection $f(x)\in L(F_{1}(x),F_{2}(x))$, there exist
continuous selections $f_{i}(x)\in F_{i}(x)$, $i=1,2$, such that
$f(x)=L(f_{1}(x),f_{2}(x))$, for all $x\in X$.
\end{Th}

\proof By the Closed graph theorem
\cite[Theorem~3.1.8]{Aubin} and \cite[Corollary 3.3.6]{Aubin},
the set-valued mapping
$$
H(x)=(F_{1}(x),F_{2}(x))\cap L^{-1}(f(x))
$$
is upper semicontinuous,
for all $x\in X$. 

Note that the mapping
which associates to each convex compact subset of $\R^{n}$ its
nearest (with respect to zero) point,
is a continuous selection of
sets in the Euclidean space, in the Hausdorff metric \cite{Daniel},
\cite{Aubin+Frankowska} and hence
$$
L^{-1}(f(x))=w(x)+\L,
$$
where $w(x)=(w_{1}(x),w_{2}(x))$ is a
continuous (at $x\in X$)
projection of zero onto $L^{-1}(f(x))$ in the 
Euclidean space
$\L_{1}\oplus \L_{2}$. 

Thus we can write
$$
H(x)=w(x)+(F_{1}(x)-w_{1}(x),F_{2}(x)-w_{2}(x))\cap\L.
$$
Consequently, we can assume that $w(x)=0$,
$H(x)=(F_{1}(x),F_{2}(x))\cap\L\ne\emptyset$ and prove the
continuity
of the last map. We shall
assume
that $H(x_{0})$ is not a
single-point, otherwise $H$ would be
continuous at the point $x=x_{0}$
for the same reason as in Example 2.\ref{Ex5}.

Suppose that $H$ is not lower semicontinuous at the point $x_{0}$.
This means that
\begin{equation}\label{L-sect}
\exists w_{0}\in H(x_{0}),\ \ep_{0}>0,\ x_{k}\to x_{0},\ \
\mbox{such that}\ w_{0}\notin H(x_{k})+B_{\ep_{0}}(0), \ \forall
k.
\end{equation}
Let $w\in\lim\sup\limits_{k\to \infty}H(x_{k})\subset H(x_{0})$;
$w\ne w_{0}$, $w,w_{0}\in\L$. Let $[u,u_{0}]\subset
F_{1}(x_{0})\cap\L_{1}$ be a projection of the segment $[w,w_{0}]$
onto $\L_{1}$ parallel to $\L_{2}$ and $[v,v_{0}]\subset
F_{2}(x_{0})\cap\L_{2}$ a projection of the segment $[w,w_{0}]$
onto $\L_{2}$ parallel to $\L_{1}$. 

By hypothesis, we have $u\ne
u_{0}$, $v\ne v_{0}$. Sets $F_{i}(x_{0})$, $i=1,2$, are strictly
convex and $\frac{u+u_{0}}{2}\in \Int F_{1}(x_{0})$,
$\frac{v+v_{0}}{2}\in \Int F_{2}(x_{0})$. Thus we can find
$\alpha>0$ such that:
$$
B^{1}_{\alpha}\left(\frac{u+u_{0}}{2}\right)\subset
F_{1}(x_{0}),\ \
B^{2}_{\alpha}\left(\frac{v+v_{0}}{2}\right)\subset F_{2}(x_{0})
$$
and
$$
B_{\alpha}\left(\frac{w+w_{0}}{2}\right)\subset
(F_{1}(x_{0}),F_{2}(x_{0})).
$$
In the last inclusion we considered the ball of norm $\max\{ \|
u\|_{\L_{1}}, \|v\|_{\L_{2}} \}$,  $(u,v)\in\L_{1}\oplus\L_{2}$.

Without loss of generality we may
assume that $\|
w-w_{0}\|>\ep_{0}$ (otherwise we reduce $\ep_{0}>0$). Let
$t=\frac{\ep_{0}}{\| w-w_{0}\|}\in (0,1)$. By a homothety with
center $w_{0}$ and the coefficient $t$ we get
that for $\tilde w =
w_{0}+\frac{t}{2}(w-w_{0})$ the
inclusion $B_{t\alpha}(\tilde
w)\subset (F_{1}(x_{0}), F_{2}(x_{0}))$ holds and $\| \tilde
w-w_{0}\|=\frac{\ep_{0}}{2}$. 

By continuity of $F_{i}$,
$i=1,2$, we get that there exists $k_{0}$ such that
$B_{\frac12 t\alpha}(\tilde w)\subset (F_{1}(x_{k}),
F_{2}(x_{k}))$, for all $k>k_{0}$ and hence we have:
$$
\tilde w \in (F_{1}(x_{k}),F_{2}(x_{k}))\cap\L\cap
B_{\ep_{0}}(w_{0}), \ \ \forall k>k_{0},
$$
i.e. $H(x_{k})\cap B_{\ep_{0}}(w_{0})\ne \emptyset$ for all
$k>k_{0}$. This contradicts (2.\ref{L-sect}).

So we have proved that $H(x)$ is continuous in the Hausdorff
metric, for all $x\in X$. Taking the Steiner point of
$H(x)=(F_{1}(x),F_{2}(x))\cap L^{-1}(f(x))$ we obtain that
$(f_{1}(x),f_{2}(x))=s\left( H(x)\right)$. \qed

\begin{Corol}\label{srtictlyConvLSC} Let $X$ be any metric space.
Let $F_{i}:X\to2^{\R^{m_{i}}}$, $i=1,2$, be ($\ep-\delta$)-lower
semicontinuous set-valued mappings with strictly convex compact
images. Let $\L_{i}=\R^{m_{i}}$, $i=1,2$ and let
$L:\L_{1}\oplus\L_{2}\to \R^{k}$ be a linear surjection such that
$\L=\ker L$ is not parallel to $\L_{i}$, $i=1,2$. Suppose that for any $x\in
X$ and  any point $f\in L(F_{1}(x),F_{2}(x))$ there exists a
pair of distinct
points $w_{1},w_{2}\in (F_{1}(x),F_{2}(x))$
such that $f=Lw_{i}$, $i=1,2$. Then for any continuous selection
$f(x)\in L(F_{1}(x),F_{2}(x))$ there exist continuous selections
$f_{i}(x)\in F_{i}(x)$, $i=1,2$, such that
$f(x)=L(f_{1}(x),f_{2}(x))$, for all $x\in X$.
\end{Corol}

\proof We can repeat word-by-word the proof from Theorem 2.\ref{srtictlyConv}
of the
lower semicontinuity
of $H$ at any point $x_{0}$.
The only difference is 
when we choose the point $w$ as an arbitrary
point
of the set $H(x_{0})\backslash\{ w_{0}\}$. Using the Michael
selection theorem \cite{MichaelII} we can choose a
continuous
selection $(f_{1}(x),f_{2}(x)) \\ \in H(x)$.\qed

Finally, we shall prove that the exact
solution of the
splitting
problem for selections takes place on the dense subset of
arguments of $G_{\delta}$-type.

\begin{Th}
 \label{Gdelta} Let $X$ be any metric space and $Y$, $Y_{i}$, $i=1,2$ any Banach
 spaces.  Let $F_{i}: X\to2^{Y_{i}}$, $i=1,2$,
be upper semicontinuous set-valued mappings with convex compact
images and suppose that the set $\cl (F_{1}(X),F_{2}(X))$ is compact.
 Let $L:Y_{1}\oplus Y_{2}\to Y$
be any
continuous linear surjection. Then for any continuous
selection $f(x)\in L(F_{1}(x),F_{2}(x))$ there exist
a
$G_{\delta}$-set $X_{f}\subset X$ and selections $f_{i}(x)\in F_{i}(x)$, $i=1,2$, continuous on the set
$X_{f}$,
such that
$f(x)=L(f_{1}(x),f_{2}(x))$, for all $x\in X_{f}$.
\end{Th}

\proof The intersection $H(x)$ of continuous mapping
$L^{-1}(f(x))$ and the
upper semicontinuous mapping
$(F_{1}(x),F_{2}(x))$ with compact images is 
upper semicontinuous
\cite{Polovinkin&Balashov,Aubin+Frankowska}. 

Moreover, the graph
$\graph
H$ is closed. By
\cite[Theorem 3.1.10]{Aubin}, 
$H(x)$ is continuous on some dense $G_{\delta}$-set $X_{f}\subset
X$. Note that $X_{f}$ is also a metric space. Applying the Michael
selection theorem \cite{MichaelII} for the
set-valued mapping
$H:X_{f}\to 2^{Y_{1}\oplus Y_{2}}$, we obtain continuous selections
$(f_{1}(x),f_{2}(x))\in H(x)$, for all $x\in X_{f}$.\qed

We conclude by some final remarks concerning $P$-sets.

\begin{Lm}
\label{lem1} Let $A\subset\R^{n}$ be any convex compact subsets and
suppose that in terms of Definition 1, for any unit vector $q$
the operator $P_{q}$ is an open map $P_{q}:A\to P_{q}A$, in the induced
topologies. Then $A$ is a
$P$-set.
\end{Lm}

\proof Suppose that $A$ is not a $P$-set. Then for some unit
vector $q$ there exists sequence $\{ x_{k}\}\subset P_{q}A$ such
that $\lim x_{k}=x_{0}$, $\lim f_{A,q}(x_{k})=\mu_{0}>f(x_{0})$.
Let $z_{0}=(x_{0},f_{A,q}(x_{0}))$, $\tilde
z_{0}=(x_{0},\mu_{0})$, $z=\frac{1}{2}(\tilde
z_{0}+z_{0})=(x_{0},\frac{1}{2}\left(\mu_{0}+f_{A,q}(x_{0}))\right)$,
$\ep=\frac{1}{4}\| \tilde
z_{0}-z_{0}\|=\frac{1}{4}\left(\mu_{0}-f_{A,q}(x_{0})\right)$. Then  $x_{k}\notin P_{q}(B_{\ep}(z)\cap A)$, for all $k$.  This
contradicts the openness of $P_{q}$.\qed

\begin{Th}
 \label{charact} Any convex compact subset $A\subset\R^{n}$ is a $P$-set
 if and only if
for any natural number $m$ and any linear map
$L:\R^{n}\to\R^{m}$, the map
$L:A\to LA$ is open in the induced topologies.
\end{Th}

\proof The openess $L:A\to LA$  was proved in~\cite{Bal}.
By Lemma 2.\ref{lem1}, we get the converse statement, it suffices
to take $L=P_{q}$.\qed

\begin{Corol}
Let $E$ be any Banach space, $\dim E=n$, and $E=L\oplus
l$, $\dim L=n-1$, $\dim l=1$. A set $A\subset E$ is a $P$-set if
and only if the projection $A$ onto $L$, parallel to $l$, is an open map
and this property holds for any pairs of subspaces $L$, $l$ with
$E=L\oplus l$, $\dim L=n-1$, $\dim l=1$.
\end{Corol}

\proof* follows from Theorem 2.\ref{charact} and the equivalence of
the Euclidean and the given norm.\qed

\bigskip

\section{Solution of approximate splitting problem for Lipschitz selections in a Hilbert space}

Let $X,Y,Y_{i}$, $i=1,2$, be infinite-dimensional Hilbert spaces.
Define $B^{X}_{\ep}(x)=\{ y\in X\ |\ \| x-y\|\le \ep\}$. Let
$L:Y_{1}\oplus Y_{2}\to Y$ be any continuous linear surjection.
We shall consider the following problem: When can
an
arbitrary Lipschitz
continuous (resp. simply Lipschitz) selection $f(x)\in
L(F_{1}(x),F_{2}(x))$ be represented in the form $f(x)=
L(f_{1}(x),f_{2}(x))$, where $f_{i}(x)\in F_{i}(x)$ are some
Lipschitz selections, $i=1,2$.

Clearly, Example 1.1 shows that there is no positive solution
of this problem in such a formulation.

We shall prove that there exists an approximate solution of the
Lipschitz
splitting problem, namely that
for any $\ep>0$, any pair of uniformly
continuous set-valued mappings $F_{i}$, $i=1,2$, with closed
convex bounded images,
and any Lipschitz selection $f(x)\in
L(F_{1}(x),F_{2}(x))$,
there exist Lipschitz selections
$f_{i}(x)\in F_{i}(x)+B^{Y_{i}}_{\ep}(0)$ such that
$f(x)=L(f_{1}(x),f_{2}(x))$.

\begin{Th}\label{unicont}
Let $(X,\varrho )$ be any
metric space and $(Y,\|\cdot\|)$ any Banach
space. Let $F_{i}:X\to 2^{Y}$, $i=1,2$, be any
set-valued mappings
with closed convex images. Let $\omega_{i}:[0,+\infty)\to
[0,+\infty)$, $i=1,2$, be the modulus of continuity for $F_{i}$,
i.e. $\lim\limits_{t\to+0}\omega_{i}(t)=0$ and
$$
\forall x_{1},x_{2}\in X\quad h(F_{i}(x_{1}),F_{i}(x_{2}))\le
\omega_{i}(\varrho (x_{1},x_{2})),\qquad i=1,2.
$$
Let $d(x)=\min\limits_{1\le i\le 2}\diam F_{i}(x)<+\infty$ for
all $x\in X$. Suppose there exist
a function $\gamma :X\to
[0,+\infty)$ and $\alpha>0$ 
such that:
\begin{equation}\label{inBall}
\gamma (x)B^{Y}_{1}(0)\subset F_{1}(x)-F_{2}(x),\qquad\forall
x\in X,
\end{equation}
\begin{equation}\label{inBallLarge}
d(x)\le \alpha \gamma (x),\qquad\forall x\in X.
\end{equation}
Then the set-valued mapping $G(x)=F_{1}(x)\cap F_{2}(x)$ is
uniformly continuous with the modulus
$\omega(t)=\max\{\omega_{1}(t),\omega_{2}(t)\}+\alpha
(\omega_{1}(t)+\omega_{2}(t))$.
\end{Th}

\proof We shall use ideas from Theorem 2.2.1 from
\cite{Polovinkin&Balashov}. Note that $G(x)\ne\emptyset$ follows
by inclusion (3.\ref{inBall}). Fix $t>1$. Choose any pair of
points $x_{1},x_{2}\in X$ and $y_{1}\in G(x_{1})$. 

We shall
prove
that there exists a point $y_{2}\in G(x_{2})$ with
\begin{equation}\label{aim1}
\| y_{1}-y_{2}\|\le t\omega (\varrho (x_{1},x_{2})).
\end{equation}
Define $\omega=\omega (\varrho (x_{1},x_{2}))$,
$\omega_{i}=\omega_{i} (\varrho (x_{1},x_{2}))$. By
the uniform continuity of $F_{i}$ it follows that:
\begin{equation}\label{uc1}
y_{1}\in F_{1}(x_{2})+t\omega_{1}B^{Y}_{1}(0),
\end{equation}
\begin{equation}\label{uc2}
y_{1}\in F_{2}(x_{2})+t\omega_{2}B^{Y}_{1}(0).
\end{equation}
Set
$d(x_{2})=\diam F_{1}(x_{2})$. By
inclusion (3.\ref{uc1}) it follows
that there exists a point $y\in
F_{1}(x_{2})$ such that
$\| y-y_{1}\|\le t\omega_{1}$.

From this and by
the inclusion (3.\ref{uc2}) we conclude that
$$
y\in F_{2}(x_{2})+t(\omega_{1}+\omega_{2})B^{Y}_{1}(0).
$$
Let
$$
\theta = \frac{\gamma (x_{2})}{\gamma
(x_{2})+t(\omega_{1}+\omega_{2})}\in [0,1).
$$
From the previous inclusion we get the inclusion
$$
\theta y\in \theta F_{2}(x_{2})+\theta
t(\omega_{1}+\omega_{2})B^{Y}_{1}(0).
$$
Keeping in mind that $\theta t(\omega_{1}+\omega_{2}) = (1-\theta
)\gamma (x_{2})$, we get
$$
(1-\theta ) \gamma (x_{2})B^{Y}_{1}(0)\subset (1-\theta
)F_{2}(x_{2})-(1-\theta )F_{1}(x_{2}),
$$
and
$$
\theta y\in \theta F_{2}(x_{2}) + (1-\theta
)F_{2}(x_{2})-(1-\theta )F_{1}(x_{2}) = F_{2}(x_{2})-(1-\theta
)F_{1}(x_{2}).
$$
Hence there exists the point $z\in F_{1}(x_{2})$ with
$$
\theta y+(1-\theta )z\in F_{2}(x_{2}).
$$
Let $y_{2}=\theta y+(1-\theta )z$. The point $y_{2}\in
F_{1}(x_{2})$ since
$y,z\in F_{1}(x_{2})$. Thus $y_{2}\in
G(x_{2})$.

From the equality $y_{1}-y_{2}=(y_{1}-y)+(1-\theta )(y-z)$ we conclude that:
\begin{equation}\label{est}
\| y_{1}-y_{2}\|\le \| y_{1}-y\|+(1-\theta )\| y-z\|\le
t\omega_{1}+(1-\theta )d(x_{2}).
\end{equation}

If $\gamma (x_{2})=0$ then $d(x_{2})=0$ and $\| y_{1}-y_{2}\|\le
t\omega_{1}$.

If $\gamma (x_{2})>0$ then from the definition of $\theta$ and
from (3.\ref{inBallLarge}) we get
$$
(1-\theta )d(x_{2})\le \frac{t(\omega_{1}+\omega_{2})}{\gamma
(x_{2})+t(\omega_{1}+\omega_{2})}\alpha \gamma(x_{2})\le t\alpha
(\omega_{1}+\omega_{2}).
$$

So by inequality (3.\ref{est}),
$$
\| y_{1}-y_{2}\|\le t(\omega_{1}+\alpha (\omega_{1}+\omega_{2})).
$$
By taking the limit $t\to 1+0$ we obtain
$$
h(G(x_{1}),G(x_{2}))\le (\omega_{1}+\alpha
(\omega_{1}+\omega_{2})).
$$
Finally, note that in the
general case we must
take
$\max\{\omega_{1},\omega_{2}\}$ instead of $\omega_{1}$, because it
may happen that $d(x_{2})=\diam F_{2}(x_{2})$.\qed

The next three propositions 3.\ref{GeomDiff}, 3.\ref{Cheb},
3.\ref{Valentine} are well known:

\begin{Pr.} \label{GeomDiff} (\cite[Lemma 5]{Ivanov+Balashov})
Let $X$ be a
Banach space and $A_{1},A_{2}\subset X$ any convex
closed bounded sets, $d=\max\{ \diam A_{1},\diam A_{2}\}$. Let
$B^{X}_{\alpha}(x_{i})\subset A_{i}$, $i=1,2$. Then for any
$\beta\in (0,\alpha)$ the following holds:
$$
h\left(A_{1}\diff B^{X}_{\beta}(0),A_{2}\diff
B^{X}_{\beta}(0)\right)\le \frac{d}{\alpha-\beta}h(A_{1},A_{2}).
$$
\end{Pr.}

\begin{Pr.} \label{Cheb} (\cite[Lemma 4]{Ivanov+Balashov}).
Let $X$ be a Hilbert space and $A_{1},A_{2}\subset X$ any convex
closed bounded sets, $A_{i}\subset B^{X}_{r}(a_{i})$, $i=1,2$.
Then $c(A_{i})\in A_{i}$, $i=1,2$, and
$$
\| c(A_{1})-c(A_{2})\|\le 2\sqrt{6rh(A_{1},A_{2})}+h(A_{1},A_{2}).
$$
\end{Pr.}

The next proposition follows by the well-known Valentine extension
theorem \cite{ValTheorem}:

\begin{Pr.} \label{Valentine} (\cite[Lemma 3]{Ivanov+Balashov})
Let $X,Y$ be Hilbert spaces and $X_{1}\subset X$ any convex subset
of $X$. Let $w:X_{1}\to Y$ be a
uniformly continuous function. Then
for any $\ep>0$, there exists a Lipschitz continuous function
$v:X_{1}\to Y$ with $\| v(x)-w(x)\|<\ep$, for all $x\in X_{1}$.
\end{Pr.}

\begin{Lm}\label{lll}
Let $X$ be a Hilbert space, $Y$ a Banach space and $L:X\to Y$ a
continuous linear surjection. Then the
operator $L$ has
the
Lip-property.
\end{Lm}

\proof Let $\ker L=\L$ and $\L^{\perp}$ be the
orthogonal subspace
of $\L$. The set-valued mapping $L^{-1}(y)$ is Lipschitz
continuous with respect to the Hausdorff distance \cite[Corollary
3.3.6]{Aubin}. Let $R(y)=\inf\limits_{x\in L^{-1}(y)}\| x\|$,
$l(y)\in L^{-1}(y)$: $\| l(y)\|=R(y)$, and
$$
G(y) = B^{X}_{2R(y)}(0)\cap L^{-1}(y).
$$
It is well known (\cite[Theorem 2.2.2]{Polovinkin&Balashov}, see
also \cite{Berdyshev, Ornelas, Aubin+Frankowska}),
that $G(y)$ is
a Lipschitz set-valued mapping with respect to the Hausdorff
distance. This also follows by
Theorem 3.\ref{unicont}.

Now consider $H(y)=G(y)\cap\L^{\perp}$. The point $l(y)$ is a
metric
projection of zero onto $L^{-1}(y)$, hence $l(y)\in G(y)$ and
(because of $l(y)\in\L^{\perp}$) $l(y)\in H(y)$. Moreover, from
the fact that some shift of $\L$ contains $G(y)$,
we can
deduce that
$H(y)=\{ l(y)\}$.

From the properties
$$
B^{X}_{\sqrt{3}R(y)}(0)\subset G(y)-\L^{\perp},\qquad \diam
G(y)\le \frac{2}{\sqrt{3}}\sqrt{3}R(y)
$$
and from Theorem 3.\ref{unicont} we now 
obtain that $H(y)=\{ l(y)\}$
is Lipschitz continuous in the
Hausdorff distance, hence $l(y)$ is a
Lipschitz function.\qed

\begin{Remark}\label{InvF}
We gave
a
purely geometric
proof of  Lemma 3.\ref{lll}. Note
that this lemma can
also be proved with the help of Implicit function theorem \cite{Ts}.
\end{Remark}

\begin{Th}\label{ApproxSolution}
Let $X,Y,Y_{i}$, $i=1,2$, be Hilbert spaces and $X_{1}\subset X$
a convex subset of $X$. Let $L:Y_{1}\oplus Y_{2}\to Y$ be a
continuous linear
surjection. Let $F_{i}:X_{1}\to 2^{Y_{i}}$, $i=1,2$, be uniformly
continuous (with modulus $\omega_{i}$) set-valued mappings with
convex closed bounded images and $d=\sup\limits_{x\in X_{1}} \diam
(F_{1}(x),F_{2}(x))<+\infty$. Suppose that for all $x\in X_{1}$ $f(x)\in
L(F_{1}(x),F_{2}(x))$ is
a Lipschitz selection. Then for any $\ep>0$
there exist Lipschitz selections $f_{i}(x)\in
F_{i}(x)+B^{Y_{i}}_{\ep}(0)$ with $f(x)=L(f_{1}(x),f_{2}(x))$, for
all $x\in X_{1}$.
\end{Th}

\proof Fix $\ep>0$. Let $f(x)\in L(F_{1}(x),F_{2}(x))$ be a
Lipschitz selection. The set-valued mapping $L^{-1}(f(x))$ is
Lipschitz continuous in Hausdorff metric \cite[Corollary
3.3.6]{Aubin}. Let $w(x) = l(f(x))\in L^{-1}(f(x))$. The function
$w(x)$ is Lipschitz continuous as a
superposition of two Lipschitz
functions: $l(y)$ (Lemma 3.\ref{lll}) and $f(x)$. 

Hence
$L^{-1}(f(x))=w(x)+\L$, $\L=\ker L$. Consider
$$
H(x) = ((F_{1}(x),F_{2}(x))+B^{Y_{1}\oplus Y_{2}}_{\ep}(0))\cap
L^{-1}(f(x))$$
$$=w(x)+((F_{1}(x),F_{2}(x))-w(x)+B^{Y_{1}\oplus
Y_{2}}_{\ep}(0))\cap\L.
$$
Without loss of generality we may assume that $w(x)=0$ and
$H(x)=((F_{1}(x),F_{2}(x))+B^{Y_{1}\oplus
Y_{2}}_{\ep}(0))\cap\L$. The mappings $x\to
(F_{1}(x),F_{2}(x))+B^{Y_{1}\oplus Y_{2}}_{\ep}(0)$  and $x\to
\L$ are uniformly continuous ($x\in X_{1}$),
%$d=\max\limits_{x\in X} \diam
%(F_{1}(x),F_{2}(x))<+\infty$,
$$
B^{Y_{1}\oplus Y_{2}}_{\ep}(0)\subset \left( (F_{1}(x),F_{2}(x))
+ B^{Y_{1}\oplus Y_{2}}_{\ep}(0)\right)-\L
$$ and
$$\diam
((F_{1}(x),F_{2}(x))+B^{Y_{1}\oplus Y_{2}}_{\ep}(0)))\le
d+2\ep\le \frac{d+2\ep}{\ep}\ep.
$$
 Using Theorem 3.\ref{unicont} we
obtain that $H(x)$, $x\in X_{1}$, is
a
uniformly continuous
set-valued mapping.

Consider $\L$ with the
induced topology: the
ball $B^{\L}_{r}(w)\subset
\L$ ($w\in \L$) is $B^{\L}_{r}(w)=B^{Y_{1}\oplus Y_{2}}_{r}(w)\cap
\L$. Obviously, $\L$ is a
Hilbert space. The set-valued mapping
$H(x)\subset \L$ has a
nonempty interior in $\L$, moreover
$H(x)\diff B^{\L}_{\ep}(0)\ne \emptyset$, for all $x\in X_{1}$.

Let $H_{1}(x)=H(x)\diff B^{\L}_{\ep/2}(0)$. By Proposition
3.\ref{GeomDiff} we have that $H_{1}(x)$ is a
uniformly continuous
set-valued mapping with convex closed images. By
Proposition
3.\ref{Cheb} we have that the Chebyshev center $c(H_{1}(x))$ of
the
set-valued mapping $H_{1}(x)$ is a
uniformly continuous function and
$c(H_{1}(x))\in H_{1}(x)$. 

By Proposition 3.\ref{Valentine} there
exists a Lipschitz continuous function $v(x)\in
c(H_{1}(x))+B^{\L}_{\ep/2}(0)$, for all $x\in X_{1}$. Hence
$$
v(x)\in c(H_{1}(x))+B^{\L}_{\ep/2}(0)\subset
H_{1}(x)+B^{\L}_{\ep/2}(0)\subset H(x).
$$
We can now
choose $(f_{1}(x),f_{2}(x))=v(x)$.\qed

\begin{Remark}\label{Ornelas}
In the finite-dimension case (when $\dim Y_{i}<\infty$, $i=1,2$)
Lip-property of $L$ follows from results \cite{Berdyshev,Ornelas}
(see also \cite{Aubin+Frankowska}, \cite{Polovinkin&Balashov},
\cite{Ivanov+Balashov}). Let $R(y)=\inf\{ \| l\|\ |\ l\in
L^{-1}(y)\}$, for all $y\in Y$. The set-valued mapping
$$
G(y) = B^{Y_{1}\oplus Y_{2}}_{2R(y)}(0)\cap L^{-1}(y)
$$
is Lipschitz continuous on $y$ (this also follows by Theorem
3.\ref{unicont}). We can choose $l(y) = s(G(y))$, where $s(\cdot)$
is the Steiner point.
\end{Remark}

\begin{Remark}\label{More}
It is easy to see that the proof of Theorem 3.\ref{ApproxSolution}
can be given in any uniformly convex Banach 
(not necessarily Hilbert)
spaces $Y_{1}, Y_{2}$,
 where every uniformly continuous function
can be approximated
(with arbitrary
precision) by a
Lipschitz function,
for any continuous linear 
surjection with Lip-property.
\end{Remark}

\begin{Example}\label{Ex7}
An exact solution of the splitting problem 
does
not exist
for Lipschitz selections. 
Besides Example 1.1 we can demonstrate
 one more example.
Tsar'kov proved \cite{Tsar'kov} that there exists
a
Lipschitz (with
respect to the Hausdorff distance) set-valued mapping $F:[0,1]\to
2^{Y}$ ($Y$ an infinite dimension Hilbert space) with convex
closed bounded images, such that the
mapping $F$ has no Lipschitz
selection. Thus for $L:Y\oplus Y\to Y$,
$L(y_{1},y_{2})=y_{1}-y_{2}$, we have $f(x)=0\in F(x)-F(x)$, but
the
Lipschitz function $f(x)=0$ cannot be represented in the form
$0=f_{1}(x)-f_{2}(x)$, where $f_{i}(x)\in F(x)$ is a Lipschitz
selection.
\end{Example}

\bigskip

\section*{Acknowledgements}
This research was supported by the
Slovenian Research Agency grants P1-0292-0101, J1-9643-0101 and
BI-RU/08-09/001. The first author was also
supported by the RFBR grant
07-01-00156.
We thank the referees for comments and suggestions.

\end{document}